\documentclass[12pt]{article}

\usepackage{amscd,amsmath, amssymb, fancyhdr, epsfig}
\usepackage[matrix,arrow]{xy}

\numberwithin{equation}{section}

\newcommand{\version}{version 0.8,\ \   01.02.2012}

\def\eqref#1{(\ref{#1})}
\newcommand{\goth}{\mathfrak}

\newcommand{\arrow}{{\:\longrightarrow\:}}
\newcommand{\Z}{{\Bbb Z}}
\newcommand{\C}{{\Bbb C}}

\newcommand{\R}{{\Bbb R}}

\newcommand{\6}{\partial}
\def\1{\sqrt{-1}\:}

\newcommand{\cntrct}                
{\hspace{2pt}\raisebox{1pt}{\text{$\lrcorner$}}\hspace{2pt}}

\makeatletter
\def\x@arrow{\DOTSB\Relbar}
\def\xlongequalsignfill@{\arrowfill@\x@arrow\Relbar\x@arrow}
\newcommand{\xlongequal}[2][]{%
        \ext@arrow 0099\xlongequalsignfill@{#1}{#2}}
\def\xlongrightarrowfill@{\arrowfill@\relbar\relbar\longrightarrow}
\newcommand{\xlongrightarrow}[2][]{%
        \ext@arrow 0099\xlongrightarrowfill@{#1}{#2}}
\makeatother

\renewcommand{\bar}{\overline}
\renewcommand{\phi}{\varphi}
\renewcommand{\epsilon}{\varepsilon}
\renewcommand{\geq}{\geqslant}
\renewcommand{\leq}{\leqslant}

\newcommand{\comass}{\operatorname{\sf comass}}

\newcommand{\Hol}{\operatorname{Hol}}

\renewcommand{\Re}{\operatorname{Re}}

\newcounter{Mycounter}[section]
\newcounter{lemma}[section]
\newcounter{claim}[section]
\newcounter{sublemma}[section]
\newcounter{corollary}[section]
\newcounter{theorem}[section]
\renewcommand{\thetheorem}{{Theorem \thesection.\arabic{theorem}}}
\newcommand{\theorem}{%
    \setcounter{theorem}{\value{Mycounter}}
    \refstepcounter{theorem}
    \stepcounter{Mycounter}
    {\noindent \bf \thetheorem:\ }}

\newcounter{conjecture}[section]
\newcounter{proposition}[section]
\renewcommand{\theproposition}
      {{Proposition \thesection.\arabic{proposition}}}
\newcommand{\proposition}{%
    \setcounter{proposition}{\value{Mycounter}}
    \refstepcounter{proposition}
    \stepcounter{Mycounter}
    {\noindent \bf \theproposition:\ }}

\newcounter{definition}[section]
\renewcommand{\thedefinition}
      {{Definition~\thesection.\arabic{definition}}}
\newcommand{\definition}{%
    \setcounter{definition}{\value{Mycounter}}
    \refstepcounter{definition}
    \stepcounter{Mycounter}
    {\noindent \bf \thedefinition:\ }}

\newcounter{example}[section]
\renewcommand{\theexample}{{Example \thesection.\arabic{example}}}
\newcommand{\example}{%
    \setcounter{example}{\value{Mycounter}}
    \refstepcounter{example}
    \stepcounter{Mycounter}
    {\noindent \bf \theexample:\ }}

\newcounter{remark}[section]
\renewcommand{\theremark}{{Remark \thesection.\arabic{remark}}}
\newcommand{\remark}{%
    \setcounter{remark}{\value{Mycounter}}
    \refstepcounter{remark}
    \stepcounter{Mycounter}
    {\noindent \bf \theremark:\ }}

\newcounter{problem}[section]
\newcounter{question}[section]
\renewcommand{\thequestion}{{Question \thesection.\arabic{question}}}
\newcommand{\question}{%
    \setcounter{question}{\value{Mycounter}}
    \refstepcounter{question}
    \stepcounter{Mycounter}
    {\noindent \bf \thequestion:\ }}

\makeatletter

\@addtoreset{equation}{section} \@addtoreset{footnote}{section} \makeatother

\def\blacksquare{\hbox{\vrule width 5pt height 5pt depth 0pt}}
\def\endproof{\blacksquare}

\begin{document}
\begin{center}
{\LARGE\bf
Subvarieties of hypercomplex manifolds with holonomy in $SL(n, {\Bbb H})$
\\[4mm]
}

Andrey Soldatenkov\footnote{Andrey Soldatenkov is partially supported by
	AG Laboratory NRU-HSE, RF government 
	grant, ag. 11.G34.31.0023}, 
Misha Verbitsky\footnote{Misha Verbitsky is partially 
supported by RFBR grant 10-01-93113-NCNIL-a,  Simons-IUM fellowship,
RFBR grant 09-01-00242-a, Science Foundation of
the SU-HSE award No. 10-09-0015 and AG Laboratory NRU-HSE, 
RF government grant, ag. 11.G34.31.0023.}

\end{center}

{\small \hspace{0.10\linewidth}
\begin{minipage}[t]{0.7\linewidth}
{\bf Abstract} A hypercomplex manifold $M$ is a manifold with
a triple $I,J,K$ of complex structure operators satisfying
quaternionic relations. For each quaternion $L=aI +bJ+cK$,
$L^2=-1$, $L$ is also a complex structure operator on $M$,
called {\bf an induced complex structure}.
We are studying compact complex subvarieties of $(M,L)$,
when $L$ is a generic induced complex structure.
Under additional assumptions (Obata holonomy 
contained in $SL(n,{\Bbb H})$, existence of an HKT metric),
we prove that $(M,L)$ contains no divisors,
and all complex subvarieties of codimension 2 are
trianalytic (that is, also hypercomplex).
\end{minipage}
}

\tableofcontents


\section{Introduction}


\subsection{Hypercomplex manifolds: an introduction}

\definition
A manifold $M$ is called {\bf hypercomplex} if $M$ 
is equipped with a triple of complex structures $I, J, K$,
satisfying the quaternionic relations 
$I\circ J =- J \circ I =K$. If, in addition, $M$ is equipped with a
Riemannian metric $g$ which is K\"ahler with respect to $I, J, K$,
it is called {\bf  hyperk\"ahler}
(\cite{_Besse:Einst_Manifo_}, \cite{_Boyer_}).

\hfill

The term ``hypercomplex manifold'' is due to C. P. Boyer,
\cite{_Boyer_}, who classified compact hypercomplex
manifolds of quaternionic dimension 1, though the
notion was considered as early as in 1955, by 
M. Obata (\cite{_Obata_}). 

The first interesting
non-hyperk\"ahler examples of hypercomplex manifolds were
found by physicists in \cite{_SSTvP_},
and independently by D. Joyce in \cite{_Joyce_}. In the same paper,
Joyce classified all homogeneous 
hypercomplex structures on simply connected
compact manifolds, using the Wang's classification
of homogeneous spaces (\cite{_Wang_}).

Next, we recall the definition of an HKT-metric. 

\hfill

\definition\label{_HKT_Definition_}
Let $(M,I,J,K)$ be a
hypercomplex manifold and $g$ a quaternionic Hermitian metric.
Consider the Hermitian forms:
$$
\omega_I(X,Y) = g(IX,Y), \quad \omega_J(X,Y) = g(JX,Y), \quad \omega_K(X,Y) = g(KX,Y).
$$
If any two of these forms are closed, the manifold is
hyperk\"ahler.  Define $\Omega_I = \omega_J + \sqrt{-1}\omega_K$.
It is easy to check that $\Omega_I \in \Lambda^{2,0}_I
M$. The metric $g$ is called HKT (``hyperk\"ahler with torsion'')
if $\6 \Omega_I = 0$, where $\6\colon \Lambda^{p,q}_I M\to \Lambda^{p+1,q}_I M$
is the $(1,0)$-part of the de Rham differential. In this case, the 
form $\Omega_I$ is called an {\bf HKT-form,}
and $(M,I,J,K,g)$ {\bf an HKT-manifold.}

\hfill

HKT-metrics were introduced by P. S. 
Howe and G. Papadopoulos \cite{_Howe_Papado_}
(see also \cite{_Gra_Poon_}) and were much studied since then.
Existence of an HKT metric puts a significant constraint
on a global geometry of a hypercomplex manifold (\cite{_Fino_Gra_}, 
\cite{_BDV:nilmanifolds_}).

Since the advent of
string theory, hypercomplex manifolds became
an important object in physics, because the
corresponding $\sigma$-models exhibit interesting
supersymmetries (\cite{_GHR_}). After A. Strominger's
paper \cite{_Strominger:Bismut_}, supersymmetric $\sigma$-models
associated with non-K\"ahler target spaces
became a popular object of study. Strominger proposed
to use the antisymmetric torsion connections
on the target spaces. In mathematics, such structures
were studied by J. Bismut \cite{_Bismut:connection_}
in connection with the local index formula.
In the hypercomplex setting, Bismut connections
were studied by P. S. Howe and G.
Papadopoulos in 1990-ies in a series of papers
starting with \cite{_Howe_Papado_}. This research 
lead them to the discovery of HKT metrics.

Since \cite{_Gra_Poon_},
HKT-metrics became an important ingredient in
the mathematical study of hypercomplex geometry.
The HKT metrics share much in common with the K\"ahler
structures. Like K\"ahler metrics, they are locally defined by a potential
(\cite{_Banos_Swann_}), but can be used to obtain Hodge-theoretic
restrictions on the cohomology (\cite{_Verbitsky:HKT_}).

In the present paper, we use the HKT-geometry 
to study complex subvarieties in hypercomplex manifolds.

\hfill

Any hypercomplex manifold admits a torsion-free connection
preserving $I, J$ and $K$, which is necessarily
unique. This connection is called {\bf the Obata
  connection}, after M. Obata, who discovered it in
\cite{_Obata_}. Any almost complex structure which is
preserved by a torsion-free connection is necessarily
integrable. Therefore, for any $a, b, c
\in \R$, with $a^2 + b^2 + c^2= 1$, the almost complex
structure $L = a I + b J + c K$ is in fact integrable. By
Newlander-Nirenberg theorem, $L$ defines a complex structure
on $M$. We denote by $(M,L)$ the complex manifold corresponding to
this complex structure.

\hfill

\definition
A complex structure $L = a I + b J + c K$, with $a^2 + b^2 + c^2= 1$,
is called {\bf induced by quaternions}, and the corresponding family,
parametrized by $\C P^1 \cong S^2$ --- {\bf the twistor family}.

\hfill

Let $(V,I,J,K)$ be a quaternionic vector space of real dimension $4n$.
The group $GL(n, {\Bbb H})$ consists of linear transformations of $V$
that preserve the complex structures $I, J$ and $K$. Consider the Hodge
decomposition $V\otimes_{\Bbb R}{\Bbb C} = V^{1,0}_I\oplus V^{0,1}_I$,
where $V^{1,0}_I$ and $V^{0,1}_I$ are eigenspaces of $I$ with eigenvalues
$\sqrt{-1}$ and $-\sqrt{-1}$ respectively. Let $\Lambda^{2n,0}_I V$
be the top exterior power of $V^{1,0}_I$. Recall that $SL(n, {\Bbb H})$
is a subgroup in $GL(n, {\Bbb H})$ consisting of those elements that
act trivially on $\Lambda^{2n,0}_I V$.

Let $(M, I, J, K)$ be a hypercomplex manifold and $\nabla$ the corresponding
Obata connection. Denote by $\Hol(\nabla)$ the holonomy group of $\nabla$.
Since the Obata connection preserves the quaternionic structure, we have
$\Hol(\nabla) \subset GL(n, {\Bbb H})$.

\hfill

\definition If the holonomy group 
$\Hol(\nabla)$ of the Obata connection on 
a hypercomplex manifold $M$ is contained in $SL(n, {\Bbb  H})$, 
we call $M$ an $SL(n, {\Bbb H})$-manifold.

\hfill

\remark
It is easy to see that an $SL(n, {\Bbb H})$-manifold
has a trivial canonical bundle (in fact, the canonical
bundle of such a manifold has a canonical flat
 connection with trivial monodromy). The converse
is also true, for compact manifolds admitting
an HKT-metric, as follows from the Hodge theory
of HKT-manifolds (\cite{_Verbitsky:canoni_}).

\hfill

\example\label{_Nil_Example_}
Let $G$ be a connected, simply connected nilpotent Lie group,
and $\Gamma \subset G$ a discrete, co-compact subgroup. 
The quotient $N := \Gamma\backslash G$
is called {\bf a nilmanifold}. Suppose that $I, J, K$ are left-invariant
complex structures on $G$ that satisfy quaternionic relations. Then the 
hypercomplex structure descends to $N$ and we call $N$ {\bf a hypercomplex nilmanifold}.
It was shown in \cite{_BDV:nilmanifolds_} that any hypercomplex nilmanifold
is in fact an $SL(n, {\Bbb H})$-manifold.

\hfill

\example
One more example of an 
$SL(n, {\Bbb H})$-manifold, due to A. Swann, is a torus fibrations
over a hyperk\"ahler base (\cite{_Swann_}). 
Let $(X,I,J,K)$ be a hyperk\"ahler manifold. A 2-form
$\alpha\in \Lambda^2 X$ is called anti-self-dual if it is of type (1,1) with
respect to any induced complex structure. If $\alpha$ represents an integral
cohomology class, then it defines a principal $U(1)$-bundle over $X$.
Given $4k$ such forms, $\alpha_1,\ldots,\alpha_{4k}$, we obtain a principal
$T^{4k}$-bundle $\pi\colon M\to X$. This bundle admits an instanton connection $A$,
given by 1-forms $\theta_i\in \Lambda^1 M$, such that $d \theta_i = \pi^*(\alpha_i)$.
The hypercomplex structure on $M$ is defined as follows: on horizontal subspaces of $A$
the quaternionic action is lifted from $X$, and on vertical subspaces it
is given by a flat hypercomplex structure of $4k$-dimensional torus.
From this construction it is easy to see that $M$ is an 
$SL(n, {\Bbb H})$-manifold. 

\hfill

The Hopf manifold 
$H = ({\Bbb  H}^n\backslash 0)/\langle A\rangle$
equipped with a standard hypercomplex structure
is not an $SL(n, {\Bbb H})$-manifold. Indeed, the
holonomy of the Obata connection on $H$ is
$\Z$ acting on $TM$ as a matrix $A$ with all eigenvalues $|\alpha_i|>1$.

It follows from the adjunction formula that
none of the homogeneous hypercomplex manifolds constructed
by Joyce in \cite{_Joyce_} has holonomy in $SL(n, {\Bbb
  H})$. Indeed, such a manifold is fibered over
a homogeneous Fano manifold with toric fibers,
hence its canonical bundle is non-trivial. However,
an $SL(n, {\Bbb  H})$-manifold has trivial
(even flat) canonical bundle.
It was shown in \cite{_Soldatenkov:SU(2)_} that
the manifold $SU(3)$ with its Joyce hypercomplex
structure has holonomy $GL(2, {\Bbb  H})$;
a similar conjecture is stated, but not proven,
for all homogeneous hypercomplex manifolds.

\subsection{Trianalytic subvarieties}

\definition
Let $M$ be a hypercomplex manifold.
A subset $Z\subset M$ is called 
{\bf trianalytic} if it is complex analytic
in $(M, L)$ for all induced complex structures
$L$.

\hfill

Geometry of trianalytic subvarieties was studied at
some length in \cite{_Verbitsky:desing_}.
It was shown that any trianalytic subvariety
can be desingularized by taking a normalization,
and this desingularization is smooth and hypercomplex.  

\hfill

The following theorem was proved in \cite{_Verbitsky:subvar_}
(see also \cite{_Verbitsky:trianalyt_}).

\hfill

\theorem\label{_complex=>triana_hk_Theorem_}
Let $M$ be a hyperk\"ahler manifold, not necessarily
compact. Then there exists a 
countable subset $S \subset \C P^1$ of induced complex structures,
such that for all compact complex subvarieties $Z \subset (M, L)$,
$L\notin S$, the subset $Z\subset M$ is trianalytic.

\hfill

\remark 
We call an induced complex structure $L$ generic,
if $L\in \C P^1 \backslash S$. If $L$ is a generic induced complex structure
on a hyperk\"ahler manifold $M$, then $(M, L)$ has no compact
complex subvarieties
except trianalytic subvarieties. Since
trianalytic subvarieties are hypercomplex
in their smooth points, their complex codimension
is even. Therefore, such $(M,L)$ has no compact
odd-dimensional subvarieties. This implies that
$(M,L)$ is not algebraic.

\hfill

It is interesting to note that this
result is manifestly false for a general
hypercomplex manifold. For example, consider
a Hopf surface $H:=({\Bbb H}\backslash 0)/(x\sim 2x)$.
For each induced complex structure $L=aI+bJ+cK$,
the manifold $(H, L)$ is fibered over $\C P^1$
with fibers elliptic curves, isomorphic to
$(\C\backslash 0)/(x\sim 2x)$. Therefore,
$(H,L)$ contains divisors for each 
induced complex structure $L$.

However, for $SL(n, {\Bbb H})$-manifolds
we still retain some control over subvarieties.
In \cite{_Gra_Verb_}, the results of
\cite{_Verbitsky:trianalyt_}
 and \cite{_Verbitsky:subvar_}
were interpreted in terms of calibrations on hyperk\"ahler
manifolds (\ref{_calibra_Definition_}).
It turns out that some of the calibrations constructed
in hyperk\"ahler geometry survive in a more general
hypercomplex setting
(\ref{_coisotro_calibra_on_HKT_Theorem_}).
This is used to obtain a weaker version
of \ref{_complex=>triana_hk_Theorem_}:

\hfill

\theorem\label{_main_Theorem_}
Let $(M, I, J, K)$ be an $SL(n, {\Bbb H})$-manifold
admitting an HKT-metric. Then there exists
a countable subset $S\subset \C P^1$, such that
for any induced complex structure $L\in \C P^1\backslash S$,
the manifold $(M,L)$ has no compact divisors, and 
all compact complex subvarieties $Z\subset (M,L)$ of complex 
codimension 2 are trianalytic.

\hfill

{\bf Proof:} See the paragraph after the 
proof of \ref{_function_cali_finite_extrema_Theorem_}. \endproof

\hfill

Without an HKT assumption, one can prove non-existence
of holomorphic Lagrangian subvarieties (for a definition
of holomorphic Lagrangian subvarieties in $SL(n, {\Bbb H})$-manifolds,
please see \ref{_Holo_Lagra_Definition_}).

\hfill

\theorem\label{_main_holo_Lag_Theorem_}
Let $(M, I, J, K)$ be an $SL(n, {\Bbb H})$-manifold. 
Then there exists a countable subset $S\subset \C P^1$, such that
for any induced complex structure $L\in \C P^1\backslash S$,
the manifold $(M,L)$ has no compact 
holomorphic Lagrangian subvarieties.

\hfill

{\bf Proof:} For any holomorphic
Lagrangian subvariety $X\subset (M,I)$, one
has $TX\cap J(TX)=0$, because $TX\subset TM$ is a
Lagrangian subspace, for any quaternionic Hermitian
metric. Therefore, \ref{_main_holo_Lag_Theorem_}
is immediately implied by 
\ref{_function_cali_finite_extrema_Theorem_} (see also \ref{_No_HKT_Remark_}). \endproof

\hfill

In the following section
we recall some facts about $SL(n, {\Bbb H})$-manifolds
and calibrations. We prove \ref{_main_Theorem_}
in Section \ref{_proofs_Section_}.


\section{Introduction to the geometry of $SL(n, {\Bbb H})$-manifolds}


This section is an introduction to HKT geometry 
of $SL(n, {\Bbb H})$-manifolds and their calibrations. We follow
\cite{_Gra_Verb_} and \cite{_Verbitsky:skoda.tex_}.

\subsection{The quaternionic
  Dolbeault complex
  on $SL(n, {\Bbb H})$-manifolds}

In this subsection, we recall the definition of a quaternionic
Dolbeault algebra of a hypercomplex manifold. 
We follow \cite{_Verbitsky:skoda.tex_}, though
this complex is essentially due to \cite{_Capria-Salamon_}.

\hfill

Let $(M,I,J,K)$ be a hypercomplex  manifold, $\dim_{\Bbb R}M=4n$.
There is a natural multiplicative action of $SU(2)\subset {\Bbb H}^*$
on $\Lambda^*(M)$, associated with the hypercomplex structure.

It is well-known that any irreducible complex representation of $SU(2)$
is a symmetric power $S^i(W_1)$, where $W_1$ is a fundamental 2-dimensional
representation. We say that a representation $U$ {\bf has
weight $i$} if it is isomorphic to $S^i(W_1)$. It follows from
the Clebsch-Gordan formula that the weight
is multiplicative in the following sense: if $i\leq j$, then
\[
W_i\otimes W_j = \bigoplus_{k=0}^i W_{i+j-2k},
\]
where $W_i=S^i(W_1)$ denotes the irreducible representation of weight $i$.

\hfill

Let $V^i\subset \Lambda^i(M)$ be a sum of all irreducible subrepresentations
$W\subset \Lambda^i(M)$ of weight $<i$. Since the weight is multiplicative,
$V^*= \bigoplus_i V^i$ is an ideal in $\Lambda^*(M)$.

It is easy to see that the de~Rham differential $d$ increases the weight
by one at most: $dV^i\subset V^{i+1}$. So $V^*\subset \Lambda^*(M)$ is
a differential ideal in the de~Rham DG-algebra $(\Lambda^*(M), d)$.

\hfill

\definition\label{_qD_Definition_}
Denote by $(\Lambda^*_+(M), d_+)$ the quotient algebra $\Lambda^*(M)/V^*$. It is called {\bf the quaternionic
Dolbeault algebra of
  $M$}, or {\bf the quaternionic Dolbeault complex}
(qD-algebra or qD-complex for short).

\hfill

The Hodge bigrading is compatible with the weight 
decomposition of $\Lambda^*(M)$, and gives a Hodge
decomposition of $\Lambda^*_+(M)$ (\cite{_Verbitsky:HKT_}):
\[
\Lambda^i_+(M) = \bigoplus_{p+q=i}\Lambda^{p,q}_{+,I}(M).
\]
The spaces $\Lambda^{p,q}_{+,I}(M)$ are the weight spaces for a particular choice of a Cartan subalgebra in
$\goth{su}(2)$. The $\goth{su}(2)$-action induces an isomorphism of the weight spaces within an irreducible
representation. This gives the following result (\cite{_Verbitsky:HKT_}):

\hfill

\proposition \label{_qD_decompo_expli_Proposition_} Let $(M,I,J,K)$ be a hypercomplex manifold and
\[
\Lambda^i_+(M) = \bigoplus_{p+q=i}\Lambda^{p,q}_{+,I}(M)
\]
the Hodge decomposition of qD-complex defined above. Then there is a natural isomorphism
\begin{equation}\label{_qD_decompo_Equation_}
	{\cal R}_{p,q}:\; \Lambda^{p+q,0}_I(M)\arrow \Lambda^{p,q}_{I,+}(M).
\end{equation}

\hfill

Consider the projection $\Pi_{+}^{p,q}\colon \Lambda^{p,q}_{I}(M)\arrow \Lambda^{p,q}_{I,+}(M)$
and let $$R:\; \Lambda^{p,q}_{I}(M)\arrow\Lambda^{p+q,0}_I(M)$$ denote the composition 
${\cal R}_{p,q}^{-1}\circ \Pi_{+}^{p,q}$.

\hfill 

Now, let $(M,I,J,K)$ be an $SL(n, {\Bbb H})$-manifold, $\dim_\R M = 4n$.
Let $\Phi_I$ be a nowhere degenerate holomorphic section of $\Lambda^{2n,0}_I(M)$.
Assume that $\Phi_I$ is real, that is, $J(\Phi_I)=\bar\Phi_I$.
Existence of such a form is equivalent 
to $\Hol(\nabla) \subset SL(n, {\Bbb H})$, where
$\nabla$ is the Obata connection (see \cite{_Verbitsky:balanced_}).
It is often convenient to define $SL(n, {\Bbb H})$-structure by
fixing the quaternionic action and the holomorphic form $\Phi_I$.

Define the map
\[ {\cal V}_{p,q}:\;
  \Lambda^{p+q,0}_I(M)\arrow\Lambda^{n+p, n+q}_I(M)
\]
by the relation
\begin{equation}\label{_V_p,q_via_test_form_Equation_}
{\cal V}_{p,q}(\eta) \wedge \alpha = \eta \wedge R(\alpha)\wedge \bar\Phi_I,
\end{equation}
for any test form $\alpha \in \Lambda^{n-p, n-q}_I(M)$.

The following proposition establishes some important properties of ${\cal V}_{p,q}$
(for the proof, see \cite{_Verbitsky:skoda.tex_}, Proposition 4.2,
or \cite{_AV:Calabi_}, Theorem 3.6):

\hfill

\proposition\label{_V_main_Proposition_}
Let $(M,I,J,K)$ be an $SL(n, {\Bbb H})$-manifold, and
\[ {\cal V}_{p,q}:\;
  \Lambda^{p+q,0}_I(M)\arrow\Lambda^{n+p, n+q}_I(M)
\]
 the map defined above.
Then
\begin{description}
\item[(i)] ${\cal V}_{p,q}(\eta)= {\cal R}_{p,q}(\eta)
  \wedge {\cal V}_{0,0}(1)$.

\item[(ii)]  The map ${\cal
V}_{p,q}$ is injective, for
  all $p$, $q$.

\item[(iii)] $(\1)^{(n-p)^2}{\cal V}_{p,p}(\eta)$ is real if and
  only $\eta\in\Lambda^{2p,0}_I(M)$ is real,
and weakly positive if and only if $\eta$ is weakly
positive.

\item[(iv)] ${\cal V}_{p,q}(\6\eta)= \6{\cal
V}_{p-1,q}(\eta)$, and ${\cal V}_{p,q}(\6_J\eta)=
  \bar\6{\cal  V}_{p,q-1}(\eta)$.

\item[(v)] ${\cal V}_{0,0}(1)= \lambda {\cal
  R}_{n,n}(\Phi_I)$, where $\lambda$ is a positive rational number,
depending only on the dimension $n$.
\end{description}


\subsection{Calibrations on $SL(n, {\Bbb H})$-manifolds}


In this subsection, we recall the construction of
the sequence of calibrations on $SL(n, {\Bbb H})$-manifolds,
following \cite{_Gra_Verb_}. These calibrations
will play the central role in the proof of the main theorem.

\hfill

\definition
(\cite{_Harvey_Lawson:Calibrated_})
Let $(V,g)$ be a Euclidean space. For
any $p$-form $\eta\in \Lambda^p (V^*)$, let $\comass(\eta)$
be  the maximum of $\frac{\eta(v_1, v_2,\ldots ,v_p)}{|v_1||v_2|\ldots |v_p|}$, for all $p$-tuples
 $(v_1, ..., v_p)$ of vectors in $V$.

\hfill

\definition\label{_calibra_Definition_}
(\cite{_Harvey_Lawson:Calibrated_})
A {\bf precalibration} on a Riemannian manifold 
is a differential form $\eta$ with $\comass(\eta) \leq 1$ everywhere.
A {\bf calibration} is a precalibration which is closed.

\hfill

Let $(M, I, J, K)$ be an $SL(n, {\Bbb H})$-manifold,
with $\Phi_I$ a holomorphic volume form on $(M,I)$
preserved by the Obata connection. We will assume that
$\Phi_I$ is real, that is $J(\Phi_I)=\bar\Phi_I$. A number of
interesting calibrations can be constructed in this situation.
The following theorem was proved in \cite{_Gra_Verb_} (Theorem 5.4):

\hfill

\theorem\label{_coisotro_calibra_on_HKT_Theorem_}
Let $(M, I, J, K)$ be an $SL(n, {\Bbb H})$-manifold, and $(\Phi_I)_J^{n,n}$ the $(n, n)$-part
of $\Phi_I$ taken with respect to $J$, and $g$ an HKT metric. Then there exists a function $c_i(m)$ on $M$,
 such that
$V_{n+i,n+i}^J:=(\Phi_I)_J^{n,n}\wedge \omega_J^i$ is
a calibration with respect to the conformal metric $\widetilde{g}=c_i g$, calibrating
complex subvarieties of $(M,J)$ which are coisotropic with
respect to the $(2,0)$-form $\widetilde{\omega}_K + \1 \widetilde{\omega}_I$.

\hfill

Note, that $V_{n+i,n+i}^J \in \Lambda^{n+i,n+i}_J M$, but using the same construction
we can obtain a similar calibration $V_{n+i,n+i}^L \in \Lambda^{n+i,n+i}_L M$
for any induced complex structure $L$.

We will need the following characterization of the form $V_{n+i,n+i}^I$
(for the proof, see \cite{_Gra_Verb_}, Remark 3.8 and Proposition 3.9):

\hfill

\proposition\label{_Calibr_Proposition_}
Let $V_{n+i,n+i}^I \in \Lambda^{n+i, n+i}(M,I)$ be a calibration
from \ref{_coisotro_calibra_on_HKT_Theorem_}. Then $V_{n+i,n+i}^I$ is
proportional to $\mathcal{V}_{i,i}(\Omega_I^i)$ and to $\Pi_{+}^{n+i,n+i}(\omega_I^{n+i})$
with some positive coefficients that do not depend on the complex structure $I$
(here $\Omega_I$ is an HKT form). 
In particular, the form $V_{n+i,n+i}^I$ is of maximal weight and
for any $\alpha \in \Lambda^{n-i,n-i}$ we have
\begin{equation}\label{_Calibr_via_test_form_Equation_}
V_{n+i,n+i}^I \wedge \alpha = a_i \Omega_I^i \wedge R(\alpha)\wedge \bar\Phi_I,
\end{equation}
where $a_i$ are some positive functions on $M$.

\hfill

\remark\label{_No_HKT_Remark_}
Note that the calibrations $V_{n+i,n+i}^I$ are constructed in the
case when the metric is HKT. However, this assumption is not
necessary for $i=0$. Since by \ref{_V_main_Proposition_}
the form ${\cal V}_{0,0}(1)$ is always closed, 
\ref{_Calibr_Proposition_} is true for $i=0$ even if the metric
is not HKT. This remark makes it possible to prove 
\ref{_main_holo_Lag_Theorem_} without the HKT assumption.

\hfill

\remark In general, the form $V_{n+i,n+i}^J$ is 
{\it not} parallel with
respect to the Obata connection.
Otherwise, since $\Phi_I$ is parallel, $\omega_J$ would also be parallel.
Then the manifold $(M,I,J,K,g)$ would necessarily be
hyperk\"ahler.
In fact, $V_{n+i,n+i}^J$ is not parallel with respect to
any torsion-free connection on $M$ (see \cite{_Gra_Poon_},
Claim 6.6).

\subsection{Holomorphic Lagrangian subvarieties in
$SL(n, {\Bbb H})$-manifolds}

Let $(M,I,J,K)$ be a $SL(n, {\Bbb H})$-manifold,
and $\Phi_J\in \Lambda^{2n,0}(M,J)$ a section
of the canonical bundle of $(M,J)$ parallel
with respect to the Obata connection. 
Since $I$ and $J$ anticommute,
$I(\Phi_J)$ is a section of $\Lambda^{0,2n}(M,J)$,
 hence satisfies $I(\Phi_J)=\alpha \bar \Phi_J$,
for $\alpha$ a complex number such that $|\alpha|=1$.
Rescaling $\Phi_J$, we can always assume that $I(\Phi_J)=\bar \Phi_J$.
Denote by $\widetilde{V}_{n,n}:=\frac {1} {n!}(\Re\Phi_J)^{n,n}_I$ the $(n,n)$-part of $\Phi_J$,
taken with respect to $I$. In \cite{_Gra_Verb_}
it was shown that $\widetilde{V}_{n,n}$ is a calibration
for any quaternionic Hermitian metric which satisfies
$|\Phi_J|=1$. The corresponding calibrated 
subvarieties were described (\cite[Proposition 5.1]{_Gra_Verb_})
as follows.

\hfill

\theorem
Let $(M,I,J,K)$ be an $SL(n, {\Bbb H})$-manifold,
$X\subset M$ a subvariety, and $\widetilde{V}_{n,n}\in \Lambda^{n,n}(M,I)$
the calibration defined above. Consider a quaternionic
Hermitian metric $h$ on $(M,I,J,K)$, and let 
$\Omega:=\omega_J+\1\omega_K$ be a (2,0)-form
constructed from $h$ as in \ref{_HKT_Definition_}.
Then the following conditions are equivalent.
\begin{description}
\item[(i)] $\widetilde{V}_{n,n}$ calibrates $X$.
\item[(i)] $X \subset (M,I)$ is a complex subvariety
which is Lagrangian with respect to $\Omega$.
\end{description}
{\bf Proof:} \cite[Proposition 5.1]{_Gra_Verb_}. \endproof

\hfill

\definition\label{_Holo_Lagra_Definition_}
Let $(M,I,J,K)$ be an $SL(n, {\Bbb H})$-manifold,
and $X\subset (M,I)$ a complex subvariety. We say that
$X$ is {\bf holomorphic Lagrangian} if it is calibrated by $\widetilde{V}_{n,n}$.

\hfill

\remark
It is remarkable that one is able to define holomorphic Lagrangian
subvarieties in the absence of a holomorphic symplectic form.
More precisely, the property of being holomorphic Lagrangian
is independent from the choice of a quaternionic Hermitian
structure which determines the $(2,0)$-form
$\Omega:=\omega_J+\1\omega_K$.


\section{Subvarieties in $SL(n, {\Bbb H})$-manifolds}
\label{_proofs_Section_}


In this subsection, we prove 
the main result of this paper (\ref{_function_cali_finite_extrema_Theorem_}),
which is used to prove \ref{_main_Theorem_}.

\hfill

Let $(M, I, J, K)$ be an $SL(n, {\Bbb H})$-manifold
equipped with an HKT-metric $g$. In the previous subsection
we have constructed a sequence of closed positive forms
$V_{n+i, n+i}^I\in \Lambda^{n+i, n+i}_I M$,
$i = 0, 1, ..., n$. We will use these forms to prove
\ref{_main_Theorem_}.

\hfill

The proof of \ref{_main_Theorem_} is based  on an
observation, which is essentially linear-algebraic.
Let $(U, I, J, K)$ be a quaternionic vector space of real dimension $4n$,
$\Phi_I \in \Lambda^{2n,0}_I (U^*)$ a complex volume form and $V_{n+i, n+i}^I$
the element of $\Lambda^{n+i,n+i}_I (U^*)$ constructed in 
\ref{_coisotro_calibra_on_HKT_Theorem_}.
Consider an $I$-invariant subspace $W\subset U$, of complex dimension
$n+i$. Note that $\dim_\C (W \cap J(W))\geq 2i$. Let 
$\xi_W\in\Lambda^{n+i,n+i}_I U$ be a volume
polyvector of $W$ (it is well defined up to a scalar multiplier).
Consider a function $\psi\colon SU(2)\to \R$ mapping $g\in SU(2)$
to $\langle V_{n+i,n+i}^I, g(\xi_W)\rangle$.
Since $W$ is $I$-invariant, $\psi$ is constant on the $U(1)$-subgroup
of $SU(2)$ associated with the complex structure $I$. This allows one
to consider $\psi$ as a function on $SU(2)/U(1)=\C P^1$.

\hfill

\proposition\label{_V_p,p_linear_algebra_Proposition_}
In the above assumptions, let $\dim_\C (W \cap J(W)) = 2k$. Then
\begin{description}
\item[(i)] If $k = i$, then $\psi$ considered as a function on $\C P^1$
has strict extremum at the point corresponding to the complex structure $I$. 
\item[(ii)] If $k > i$, then $\psi$ is identically zero.
\end{description}

{\bf Proof:} Let us fix a quaternionic-hermitian metric in $U$,
such that $\Phi_I \wedge \bar{\Phi_I}$ is its volume form.
Denote by $\eta_W \in \Lambda^{n-i,n-i}(U^*)$ the form dual to $\xi_W$,
that is $\eta_W = *(\xi_W^\sharp)$, where $\sharp$ denotes the duality
with respect to the metric and $*$ is the Hodge star operator.

Then we have $\psi(g) = \langle V_{n+i,n+i}^I, g(\xi_W) \rangle = *(V_{n+i,n+i}^I\wedge g(\eta_W))$
and by (\ref{_Calibr_via_test_form_Equation_}) we obtain
$$V_{n+i,n+i}^I \wedge g(\eta_W) = a_i\Omega_I^i\wedge R(g(\eta_W))\wedge \bar{\Phi_I}.$$

We can choose an orthonormal basis in $U^{1,0}$ of the form 
\[ \langle e_1, J\bar{e_1}, \ldots , e_n, J\bar{e_n} \rangle,\]
such that 
$$W^{1,0} = \langle e_1, J\bar{e_1}, \ldots , e_k, J\bar{e_k}, e_{k+1}, e_{k+2}, \ldots, e_{n+i-k}\rangle,$$

If $k>i$ then for any $g\in SU(2)$ we see that $R(g(\eta_W))$ has to belong to the $(2n-2i)$-th
exterior power of the subspace in $(U^*)^{1,0}$ spanned by $e_{k+1}^*, J\bar{e_{k+1}^*}, \ldots, e_{n}^*, J\bar{e_{n}^*}$.
But this exterior power vanishes, so $R(g(\eta_W)) = 0$, which proves the second part of the proposition.

If $k=i$ then 
$$\eta_W = J\bar{e_{i+1}^*}\wedge \ldots \wedge J\bar{e_{n}^*}\wedge J{e_{i+1}^*}\wedge \ldots \wedge J{e_{n}^*}$$
and $$R(\eta_W) = e_{i+1}^* \wedge J\bar{e_{i+1}^*} \wedge \ldots \wedge e_{n}^* \wedge J\bar{e_{n}^*}.$$
Since $\Omega_I = \sum e_j^* \wedge J\bar{e_j^*}$, we see that in this case $\Omega_I^i\wedge R(\eta_W)$
does not vanish, that is $\psi(1) \neq 0$. We claim that the function $\psi$ is non-constant in this case.
Otherwise, the function would be equal to its average over $SU(2)$, which equals 
$\langle {\mathrm Av}_{SU(2)}V_{n+i,n+i}^I, \xi_W \rangle$. But in the last expression 
${\mathrm Av}_{SU(2)}V_{n+i,n+i}^I = 0$ because $V_{n+i,n+i}^I$ is of the maximal weight and
belongs to non-trivial irreducible representation of $SU(2)$. Therefore, the average
of $\psi$ is zero.

Now, consider the action of $U(1)$ on the two-dimensional sphere by rotations around
the axis that passes through the two points corresponding to the complex structures $I$ and $-I$.
We claim that the function $\psi$ is invariant under this action. This follows
from the definition of $\psi$: observe that $V_{n+i,n+i}^I$ and $\xi_W$ are
$I$-invariant. Therefore, $g\mapsto \langle V_{n+i,n+i}^I, g(\xi_W)\rangle$
is invariant (as a function on $SU(2)$) under the adjoint action of
the $U(1)$-subgroup corresponding to~$I$.

Since $\psi$ is an analytic non-constant 
function on the sphere, and it is invariant
under the $U(1)$-action considered above, it must have strict extremum at $I$.
This proves the proposition.
\endproof

\hfill

Let $(M, I, J, K)$ be an $SL(n, {\Bbb H})$-manifold
equipped with an HKT-metric, 
and $[Z]\in H_{2n+2i}(M,\Z)$ an integer homology class.
Consider a function 
$\phi_Z:\; \C P^1\arrow \R$ associating to each $L \in \C P^1$ 
a number $\int_Z V^L_{n+i,n+i}$, 
where $V^L_{n+i,n+i}\in \Lambda^{n+i, n+i}(M,L)$
is the corresponding calibration form.

Note that if $L = aI+bJ+cK$ with $a^2+b^2+c^2 = 1$, then 
$\omega_L = a\omega_I +b\omega_J + c\omega_K$. Since $V^L_{n+i,n+i}$ is
proportional to $\Pi_{+}^{n+i,n+i}(\omega_L^{n+i})$ with the coefficient
that does not depend on $L$ (see \ref{_Calibr_Proposition_}),
we see that the function $\phi_Z$ is a restriction to $S^2$ of a homogeneous polynomial
in $\Bbb R^3$. Such a function can have only a finite number of strict extrema
(this follows from the fact that real algebraic variety can
have only finitely many connected components, see \cite{_Whitney_}).
Let $S\subset \C P^1$ be the set of all strict extrema of $\phi_Z$ for all integer homology classes.
Since for each fixed $[Z]\in H_{2n+2i}(M,\Z)$ the set of strict extrema of $\phi_Z$ is
finite, the set $S$ is countable.

\hfill

\theorem\label{_function_cali_finite_extrema_Theorem_}
For each $L\in \C P^1\backslash S$,
and any compact complex subvariety $Z\subset (M,L)$ of complex dimension $n+i$,
one has $\dim_\C (TZ \cap J(TZ))>2i$.

{\bf Proof:} Fix a volume form $dv$ on $Z$ and assume that $\dim_\C (TZ \cap J(TZ))=2i$.
Then 
$$\phi_Z(L_1) = \int_Z \langle V^{L_1}_{n+i,n+i}, \xi_{TZ}\rangle dv.$$
Fix an arbitrary smooth point $x\in Z$. Note that for any $g\in SU(2)$ we have
$\langle V^{L}_{n+i,n+i}, g(\xi_{T_x Z})\rangle = \langle V^{{\mathrm Ad}_g L}_{n+i,n+i}, \xi_{T_x Z}\rangle$. 
Thus, \ref{_V_p,p_linear_algebra_Proposition_} implies that
the function $L_1\mapsto \langle V^{L_1}_{n+i,n+i}, \xi_{T_x Z}\rangle$
has strict extremum at $L_1 = L$. Thus, $\phi_Z$ also has strict extremum at
this point, which contradicts our assumption that $L\in \C P^1\backslash S$.
\endproof

\hfill

{\bf Proof of \ref{_main_Theorem_}:} Let $L\in \C P^1\backslash S$ and $Z$
be a divisor in $(M, L)$, that is a compact $L$-complex subvariety of
complex dimension $2n-1$. Then \ref{_function_cali_finite_extrema_Theorem_}
implies that $\dim_\C (TZ \cap J(TZ))>2n-2$. Since $TZ \cap J(TZ)$ is ${\Bbb H}$-invariant,
the last inequality would imply that the dimension equals $2n$, which is impossible.
So there exist no divisors in $(M, L)$.

Analogously, if $\dim_\C Z = 2n - 2$, then we have $\dim_\C (TZ \cap J(TZ))>2n-4$.
This implies that $\dim_\C (TZ \cap J(TZ))=2n-2$, that is, $TZ$ is ${\Bbb H}$-invariant
and $Z$ trianalytic. This completes the proof of the theorem.
\endproof

\hfill

\remark\label{_conject_Remark_}
We should note that the existence of an HKT-metric was 
essential for the proof of the main theorem. It still remains
unclear if this condition could be removed.

On the other hand, the condition that the holonomy of the Obata
connection is contained in $SL(n, {\Bbb H})$ is known to be necessary.
There exist examples of HKT-manifolds with odd-dimensional
complex subvarieties for each induced complex structure. An HKT-structure
on compact Lie groups due to D. Joyce (\cite{_Joyce_}) gives such an example: 
it is well-known (see e.g. \cite{_Verbitsky:toric_}) that these manifolds admit
a toric fibration over a rational base, hence they always
contain divisors.

\hfill

\remark
Let $T$ be a compact hyperk\"ahler torus,
and $L$ a generic induced complex structure. Then
all complex subvarieties of $T$ are again tori
(\cite{_KV:GK_}).
We conjecture that something similar would
happen for nilmanifolds, and for flat hypercomplex
manifolds.

\hfill

\question
Let $M$ be a compact $SL(n, {\Bbb H})$-manifold with flat Obata
connection, and $L$ a generic induced complex structure.
Is it true that all the complex 
subvarieties of $(M,L)$ are also flat? 

\hfill

\question
Let $M$ be 
a hypercomplex nilmanifold  (\ref{_Nil_Example_}),
and $L$ a generic induced complex structure. 
Is it true that all complex subvarieties of $(M,L)$
are also nilmanifolds?

\hfill

\hfill

\hfill

\small{

\noindent
{\sc Andrey Soldatenkov}\\
{\sc Laboratory of Algebraic Geometry,\\ 
National Research University Higher School of Economics,\\ 
7 Vavilova Str., Moscow, Russia, 117312}

\hfill

\noindent
{\sc Misha Verbitsky}\\
{\sc Laboratory of Algebraic Geometry,\\ 
National Research University Higher School of Economics,\\ 
7 Vavilova Str., Moscow, Russia, 117312}\\
{\tt verbit@maths.gla.ac.uk, \ \  verbit@mccme.ru\\}

 }

\end{document}